\documentclass{article} \font \Bbbten=msbm10 \font \Bbbsev=msbm7
 \font \Bbbfiv=msbm5
\newfam\Bbbfam \def \Bbb {\fam\Bbbfam\Bbbten}\textfont\Bbbfam =\Bbbten
\scriptfont\Bbbfam=\Bbbsev \scriptscriptfont\Bbbfam=\Bbbfiv

\usepackage{graphicx}
\usepackage{psfrag}

\begin{document}

%Teoremas-proposiciones-lemas-observaciones

\newcounter{theorem}[section]
\newtheorem{defi}[theorem]{\sc Definition}
\newtheorem{lema}[theorem]{\sc Lemma}
\newtheorem{prop}[theorem]{\sc Proposition}
\newtheorem{cor}[theorem]{\sc Corrollary}
\newtheorem{teo}[theorem]{\sc Theorem}
\newtheorem{obs}[theorem]{\sc Remark}
\def\bp{\noindent{\it Proof. }}
\def\ep{\noindent{\hfill $\fbox{\,}$}\bigskip\newline}
\renewcommand{\theequation}{\arabic{section}.\arabic{equation}}
\renewcommand{\thetheorem}{\arabic{section}.\arabic{theorem}}
%Comandos-palabras-definiciones
\newcommand{\eps}{\varepsilon}
\newcommand{\disp}[1]{\displaystyle{\mathstrut#1}}
\newcommand{\fra}[2]{\displaystyle\frac{\mathstrut#1}{\mathstrut#2}}

\title{There are no stable points for continuum-wise expansive homeomorphisms}
\author{Jana Rodriguez Hertz}
\date{}
\maketitle

%Font para hacer los Reales
\begin{abstract}
We obtain some results about continuum-expansive expansive homeomorphisms, such as
non-existence of stable points and presence of non-trivial
connected components within the local stable and unstable sets. These facts have been of importance in
theorems of classification of expansive homeomorphisms (see below).
They are achieved now, however, by new and self contained techniques that hold on more general metric spaces.
\end{abstract}
%______________________________________________________________________________________________
\section{Introduction}
%______________________________________________________________________________________________
A homeomorphism $f$ on a compact metric space $(X,d)$ is said to
be {\em expansive} if there exists a positive constant $\alpha>0$,
called an {\em expansivity constant for $f$}, such that
$\sup\{d(f^n(x),f^n(y)): n\in {\Bbb Z}\}>\alpha$ for all $x\ne y$
in $X$. This class of homeomorphisms includes, among other
examples, subshifts of finite type, pseudo-Anosov homeomorphisms
on surfaces and Anosov diffeomorphisms. Expansiveness
is invariant under topological conjugation.\par
One may easily
construct trivial examples of expansive dynamics, such as the
restriction of the North Pole-South Pole diffeomorphism to the
closure of a single (non trivial) orbit. However, the presence of
some topological properties on the space $X$ imposes severe
restrictions to the behavior of an expansive homeomorphism $f$,
and viceversa.\par
For instance, if $X=M^n$ is an $n$-dimensional orientable compact
closed manifold, then the fact that $f$ is expansive on $M$
implies:
\begin{enumerate}
\item $n\not=1$.
\item If $n=2$ then either $M={\Bbb T}^2$, in which case $f$ is
conjugated to an Anosov diffeomorphism, or else $g(M)>1$ and $f$
is then conjugated to a pseudo-Anosov diffeomorphism \cite{h,l2}.
\item If $n=3$ and $f\in C^{1+\eps}$, then $\Omega(f)=M$ implies
$f$ is conjugated to an Anosov diffeomorphism and $M={\Bbb T}^3$
\cite{v}.
\item $M$ cannot be simply connected \cite{m}.
\end{enumerate}
Results in \cite{h,l2,v} make strong use of the fact that such
homeomorphisms admit no stable points and that, furthermore, there
are non-trivial, uniformly sized, connected pieces within both the
local stable and unstable sets of each point (definitions are in
next paragraph).\par
On \cite{ka}, a weaker notion of expansiveness is introduced,
which is the following: a homeomorphism $f:X\rightarrow X$ is said
to be {\em continuum-wise expansive (cw-expansive)} if there
exists a constant $\alpha>0$ such that each non trivial continuum
$C$ satisfies $\sup_{n\in{\Bbb Z}}{\rm diam}\:
f^n(C)>\alpha$.\par
This class contains the class of expansive homeomorphisms, but it
is strictly larger, as may be seen in the following example:\par
Considering ${\Bbb S}^2$ as the quotient space of the projection
map $\pi:{\Bbb T}^2\to{\Bbb S}^2$ defined for each $x\in{\Bbb T}^2$ by
$\pi^{-1}(\pi(x))=\{x,-x\}$, and
providing ${\Bbb S}^2$ with the metric $d_{{\Bbb S}^2}(\xi,
\eta)=\min\{d_{{\Bbb T}^2}(x,y): \pi(x)=\xi\,,\pi(y)=\eta\}$;
it is not difficult to show that ${\rm diam}_{{\Bbb
S}^2}\:\pi(C)\geq \frac{1}{2}{\rm diam}_{{\Bbb T}^2}\: C $ for any
continuum $C\subset {\Bbb T}^2$ such that ${\rm diam}_{{\Bbb
T}^2}\: C<\frac 12$. Now, if $f$ is a hyperbolic
automorphism defined on ${\Bbb T}^2$, then $f$ is cw-expansive. The previous comment
implies that $F=\pi\circ f\circ\pi^{-1}$ is
a cw-expansive homeomorphism on ${\Bbb S}^2$ which cannot
be expansive due to \cite{h,l2}.\newline\par
The purpose of this paper is to provide self contained proofs
of the presence of non trivial connected pieces within the local
stable and unstable sets of each point for cw-expansive homeomorphisms acting on
locally connected compact metric spaces. The same result was obtained in \cite{h,l2}
for expansive homeomorphisms, though we use less sophisticated
machinery. As a particular case, it follows that there are no stable points
for this kind of dynamics, a property which is also known in
the literature as {\em sensitive dependence on the initial conditions}.
The incidence of different forms of connectedness in these results is
studied, obtaining counterexamples in some cases.\par%
These new tools leave the possibility
of improving classification results for this class of
homeomorphisms.
\newline\par {\em Acknowledgements:} I would like to thank people
of the Seminar of Dynamical Systems at the IMERL, particularly J.
Lewowicz for useful comments. I wish to thank the Laboratoire of
Topologie at Dijon, where part of this work was done, and
specially C.Bonatti, JM Gambaudo and E. P\'{e}cou, for kind
hospitality during my stay. K. Hiraide has gently informed me
about existing results.

%________________________________________________________________

\subsection{Main results}\label{main results}
%________________________________________________________________
If $f$ is a homeomorphism on a compact metric space,
the {\em $\eps$-local
stable set of a point $x$ in $X$} is defined as the set
$$W^s_\eps(x)=\{y\in X:\: d(f^n(x), f^n(y))\leq\eps\quad\forall n\geq 0\}$$

We shall say that $x\in X$ is a {\em stable point (in the future)}
if $\{W^s_\eps(x)\}_{\eps>0}$ is a neighborhood basis for $x$,
that is, for each $\eps>0$ there exits $\delta>0$ such that
$\sup_{n\to\infty} d(f^n(x),f^n(y))\leq\eps$ if $d(x,y)\leq
\delta$. The notions of {\em $\eps$-local unstable set} and {\em
stable point in the past} are defined similarly. Denoting by $X'$
the set of accumulation points of $X$, we show
\begin{teo}\label{no.estables.cw} If $f$ is a cw-expansive homeomorphism on
a locally connected compact metric space $X$, then there are no
stable points for $f$ in $X'$.
\end{teo}
\par Previous theorem is used to show a much stronger fact, namely,
that each point in $X'$ belongs to a non-trivial connected
component of its own stable and unstable set. This is an essential
step in \cite{h,l2,v}.

In fact, denoting by $CW^\sigma_\eps(x)$ the connected component of
$x$ in the set $W^\sigma_\eps(x)$ with $\sigma=s,u$, we show:
\begin{teo} \label{piezas.conexas.cw} If $f$ is a cw-expansive homeomorphism on a
locally connected compact metric space $X$, then for each $\eps>0$
there exists $\delta>0$ such that $$\inf_{x\in X'}{\rm
diam}\:CW^s_\eps(x)\geq\delta\qquad\mbox{ and }\qquad\inf_{x\in
X'}{\rm diam}\:CW^u_\eps(x)\geq\delta$$\end{teo}
Let us remark that in \cite{rr} there is an example showing that Theorems
\ref{no.estables.cw} and \ref{piezas.conexas.cw} do not hold on general
continua. Moreover, Theorem \ref{piezas.conexas.cw} does not hold on general
continua even under the assumption of non existence of stable points,.
An example showing this is built in \S \ref{ejemplo.1}.
%______________________________________________________________________________________________
\section{Connectedness and stability}\label{conexion+estabilidad}
%______________________________________________________________________________________________
%______________________________________________________________________________________________
\subsection{Expansive homeomorphisms}
%______________________________________________________________________________________________

In this paragraph we shall prove that there are no stable accumulation points for an
{\em expansive} homeomorphism acting on a locally connected metric space $X$ (Theorem
\ref{no.estables.exp}). Many of the steps involved, however, are valid under
much weaker hypotheses on $X$. Namely, we shall see
\begin{enumerate}\label{sketch de dem conexion+estabilidad}
\item Any stable recurrent point is periodic.
\item The set of stable points is open.
\item Any stable point having a point of local connectedness in
its $\alpha$-limit point is recurrent, and hence periodic.
\end{enumerate}
These statements are valid on {\em any compact metric space}. The
three of them together, in a locally connected setting, imply non
existence of stable accumulation points. The following property is
basic in the whole proof.
\begin{lema}\label{N.eps}$f$ is {\em uniformly expansive}, i.e.
for each $\eps>0$ there exists a positive integer $N_\eps$ such
that for either $x,y$ in $X$
$$d(x,y)\geq\eps\qquad \Rightarrow
\qquad \sup_{|n|\leq N_\eps}d(f^n(x),f^n(y))>\alpha $$
\end{lema}
\bp Otherwise, one could choose $\eps>0$ and pairs of points
$x_j,y_j\in X$ with $d(x_j,y_j)\geq\eps$, so that
 $d(f^n(x_j),f^n(y_j))\leq\alpha$ for all $|n|\leq j$, and $j\in{\Bbb Z}^+$.
Assuming that $x_j\to x_*$ and $y_j\to y_*$, we would obtain that
$\eps\leq\sup_{n\in {\Bbb Z}}d(f^n(x_*),f^n(y_*))\leq \alpha$,
contradicting expansivness. \ep
An immediate, though interesting, corollary of this is that
$\alpha$-local stable sets are {\em uniformly} shrunk under the
action of $f$. Indeed, if $y\in W^s_\alpha(x)$ then, for each
$n\geq N_\eps$ and $|j|\leq N_\eps$ we have
$d(f^{n+j}(x),f^{n+j}(y))\leq\alpha$. Lemma \ref{N.eps} implies
$d(f^n(x),f^n(y))<\eps$. This proves:
\begin{cor}\label{achicamiento.estables} For all $\eps>0$ and $x\in
X$,  $\quad f^n(W^s_\alpha(x))\subset B_\eps(f^n(x))\quad$  if
$n\geq N_\eps$.
\end{cor}
 In particular, $\omega(y)=\omega(x)$ for all
$y\in W^s_\alpha(x)$, where $\omega(x)$ denotes, as usual, the set
of $\omega$-limit points of the orbit of $x$.
We may obtain from this fact the following relationship between
stability and recurrence. Recall that $x$ is said to be a {\em
recurrent point} if $x\in\omega(x)$.
\par
\begin{prop}\label{PER}
A stable recurrent point is always periodic.
\end{prop}
\bp It suffices to consider $\delta>0$ so
that $B_\delta(x)\subset W^s_\alpha(x)$, and $m\geq N_{\delta/4}$ so that $f^m(x)\in
B_{\delta/4}(x)$. Then $f^m(B_\delta(x))\subset
f^m(W^s_\alpha(x))\subset B_{\delta/4}(f^m(x))\subset
B_{\delta/2}(x)$.\par This implies the existence of a periodic point
$z=\bigcap_{k\geq 0}\overline{ f^{km}(B_\delta(x))}$ in
$B_{\delta/2}(x)$, which means $x\in\omega(x)=\omega(z)=o(z)$ is
actually a periodic point.\ep
For the sake of simplicity, we shall denote \vspace*{-0.4em}
\begin{equation}\label{conj.estable.finito}
W_{\eps,N}^s(x)=\bigcap_{n=0}^N\{y\in X: d(f^n(x),f^n(y))\leq\eps\}
\end{equation}\vspace*{-0.4em}
We have the following characterizations of stable points.
\begin{lema}\label{equiv.estables} Each condition below is
equivalent to $x$ being a stable point
\begin{enumerate}
\vspace*{-0.4em}
\item \label{alfa.set} $W^s_\alpha(x)$ is a neighborhood of $x$
\vspace*{-0.4em}
\item \label{delta.set} There exists $0<\delta_0<\alpha$ such that for each $0<\delta\leq
\delta_0$ we have
$W^s_{\delta,N_\delta}(x)=W^s_\delta (x)$.
\end{enumerate}
\end{lema}
\bp It follows by taking $\delta_0>0$ so that $B_{\delta_0}(x)\subset
W^s_\alpha(x)$ and applying Corollary
\ref{achicamiento.estables}.\ep
As a result, if $x$ is a stable point and
$0<\delta<\min\{\delta_0,\alpha/2\}$, then all points $y$ in
${\rm int}(W^s_\delta (x))$ satisfy $W_\alpha^s(y)\supset
{\rm int} (W^s_\delta (x))$. Item \ref{alfa.set}. in
proposition above implies ${\rm int}(W_\delta^s(x))$ consists of
stable points. Hence,
\begin{cor}\label{estable.abierto}
The set of stable points is open.
\end{cor}
We shall impose now more restrictions on the set
$X$. The following property is key in showing Theorem
\ref{no.estables.exp}. We shall denote by $CB_\eps(x)$ the
component of $x$ in the ball $B_\eps(x)$
\begin{prop}\label{PECL}
Let $x$ be a stable point. If $\alpha(x)$ contains a point of
local connectedness of $X$, then $x\in\omega(x)$, and hence it is periodic.
\end{prop}
\bp The proof is reduced to proving
there exists $\rho>0$ such that:
$$B_\rho(z)\subset W_\alpha^s(f^{-n_k}(x)) \qquad\forall
\,k\geq 0$$
where $z=\lim_{k\to\infty}f^{-n_k}(x)$ is a point of local connectedness of
$X$. Indeed, this would give us for each $\eps>0$, a pair
$n_k>m_k\geq
N_\eps$ such that $f^{-m_k}(x)\in B_\rho(z)\subset
W_\alpha^s(f^{-n_k}(x))$. Using Corollary
\ref{achicamiento.estables} we would immediately have that $f^{n_k-m_k}(x)\in
B_\eps(x)$ whence
$x\in\omega(x)$.\par
If a $\rho>0$ as described above did not exist, we would find a
subsequence $f^{-n_k}(x)\in CB_{1/k}(z)$ such that
$CB_{1/k}(z)\not\subset W_\alpha^s(f^{-n_k}(x))$. Choosing
$0<\delta<\alpha$ as in Lemma \ref{equiv.estables}, so that
$W^s_\delta(x)=W^s_{\delta, N}(x)$, we would find, via
connectedness, a point $f^{-n_k}(y_k)$ in $CB_{1/k}(z)\cap\partial
W_\delta^s(f^{-n_k}(x))$, that would satisfy
$\sup_{n\geq-n_k}d(f^n(x),f^n(y_k))=d(f^{m_k-n_k}(x),f^{m_k-n_k}(y_k))=\delta$
for some $m_k\geq 0$ verifying $m_k\to+\infty$, due to continuity
of $f$.\par We may assume that $f^{m_k-n_k}(x)$,
$f^{m_k-n_k}(y_k)$ converge, respectively, to $x^*$, $y^*$, whence
it would follow that $d(x^*,y^*)=\delta$. On the other hand, for
each $l\in {\Bbb Z}$, we would have
$$d(f^l(x^*),f^l(y^*))=\lim_{k\to\infty}
d(f^{l+m_k-n_k}(x),f^{l+m_k-n_k}(y_k))\leq \sup_{n\geq
-n_k}\{d(f^n(x),f^n(y_k))\}=\delta$$  contradicting
the hypothesis of expansiveness.\ep
As an immediate corollary we obtain:
\begin{teo}\label{no.estables.exp}
If $f$ is an expansive homeomorphism over a locally connected
compact metric space $X$, then there are no stable points for $f$
in $X'$.
\end{teo}
We observe that, though it is beyond the scope of this paper,
Proposition \ref{PECL} could also be useful in identifying
``basins of attraction" in general continua. For instance, in such
a basin, all but possibly one point shall have $\alpha$-limit sets
consisting of non locally connected points.\newline\par
%

%_____________________________________________________________________
\subsection{Continuum-wise expansive homeomorphisms}
%______________________________________________________________________
%
If $f:X\to X$ is an $\alpha$ cw-expansive
homeomorphism on a compact metric space $X$, then a
little make-up allows us to extend Theorems \ref{no.estables.exp}
to the cw-expansive case. Essentially the same steps will be followed, though
special care is required at some steps.
\begin{lema}\label{cw-uniform} $f$ is {\em uniformly cw-expansive}, that is,
for each $\eps>0$ there exists a positive integer $N_\eps$ such
that $\max_{|n|\leq N_\eps}{\rm diam}\:f^n(C)>\alpha $ for any continuum $C$
satisfying ${\rm diam}\:C\geq\eps$.
\end{lema}
The proof follows substantially as in Lemma \ref{N.eps}, as well as
the following corollary:
\begin{cor}\label{achicamiento.estables.cw} For all $\eps>0$ and $x\in
X$, $f^n(CW^s_{\alpha/2}(x))\subset B_\eps(f^n(x))$ if $n\geq
N_\eps$. In particular, $\omega(y)=\omega(x)$ for all $y\in
CW^s_{\alpha/2}(x)$.
\end{cor}
As in \S2, we shall prove items 1.,2. and 3. stated in page
\pageref{sketch de dem conexion+estabilidad}; however, we must
take into account that, unlike in the expansive case, locally
connectedness is required. Observe that if $X$ is locally
connected at $x$, then $x$ being a stable point is equivalent to the fact
that $\{CW_\eps^s(x)\}_{\eps>0}$ is a neighborhood basis of $x$.
Having this in mind, a procedure very similar to that in proof of
Proposition \ref{PER} yields
\begin{prop} Any
stable recurrent point of local connectedness of $X$ is periodic.\end{prop}
Observe that if $X$ is locally connected at $x$, then all $y\in{\rm
int}\,CW_{\alpha/2}(x)$ are stable. Indeed, for each $\eps>0$ there
exists $\rho>0$ satisfying $B_\rho(y)\subset {\rm
int}\,CW_{\alpha/2}(x)\cap W_{\eps, N_{\eps/2}}^s(y)$. Corollary
\ref{achicamiento.estables.cw} implies $\sup_{n\geq 0}d(f^n(y),
f^n(z))\leq\eps$. This proves:
\begin{prop} \label{estable.abierto.cw} If $X$ is locally connected at a stable point $x$, then
there exists a neighborhood of $x$ consisting of stable points. In
particular, the set of stable points is open if $X$ is locally
connected.
\end{prop}
The following characterization of stable points shall also be used
\begin{lema}\label{equivalencias.cw} If $X$ is locally connected at $x$, each condition below is
equivalent to the stability of the point
$x$.
\begin{enumerate}
\item \label{alfa.set.cw} $CW^s_{\alpha/2}(x)$ is a neighborhood of
$x$ \vspace*{-0.4em} \item \label{delta.set.cw} There exists
$0<\delta_0<\alpha$ such that for all $0<\delta\leq \delta_0$,
$W_{\delta, N_\delta}^s(x)=W^s_\delta (x)$
\end{enumerate}
\end{lema}
The following proposition directly implies Theorem \ref{no.estables.cw}.
\begin{prop}\label{estable.periodico.cw}
If $X$ is locally connected, then any stable point $x$ is
recurrent and hence, periodic.
\end{prop}
\bp Following the spirit of the proof of Proposition \ref{PECL}, it will suffice to
see that there is some $\rho>0$ satisfying
$B_\rho(f^{-n}(x))\subset CW^s_{\alpha/2}(f^{-n}(x))\qquad \mbox{for all }n\geq
0$.\par
If it were not the case, we
would find a subsequence $n_k\to\infty$ such that $CB_{1/k}(f^{-n_k}(x))\setminus CW_{\alpha/2}(f^{-n_k}(x))$
is not empty. Choosing $\delta>0$ as in Lemma \ref{equivalencias.cw} so that $W_\delta^s(x)\subset
CW_{\alpha/2}^s(x)$ we would find, due to connectedness of
$CB_{1/k}(f^{-n_k})$, a point $f^{-n_k}(y_k)\in
CB_{1/k}(f^{-n_k}(x))$ belonging to $\partial
CW_\delta(f^{-n_k}(x))$. Local connectedness of $X$ would imply that
$f^{-n_k}(y_k)$ actually belongs to $\partial
W_\delta(f^{-n_k}(x))$, and hence
$d(f^{m_k-n_k}(x),f^{m_k-n_k}(y_k)=\delta$
for some $m_k\geq 0$, which, because of continuity, would satisfy $m_k\to \infty$.
This would mean that
$C_k=f^{m_k}(CW_\delta(f^{-n_k})(x))$ is a continuum satisfying
$$\delta\leq \sup_{n\geq -m_k}{\rm diam } f^n(C_k)\leq 2\delta\leq \alpha$$
Assuming that $C_k$ converges in the Hausdorff topology to a
continuum $C_*$, we would obtain that $\delta\leq \sup_{n\in {\Bbb
Z}}{\rm diam } f^n(C_*)\leq \alpha$, what would mean that $f$ is
not cw-expansive.\ep Now we turn our attention to proving Theorem
\ref{piezas.conexas.cw}.  The following Lemma is straightforward:
\begin{lema}\label{estabilidad.finita}If $X$ is a locally connected compact metric space,
then for all $0<\eps<\alpha$ and $\delta>0$ there exists
$K=K_{\eps,\delta}>0$ such that all $x\in X'$ satisfy
$$CB_\delta(x)\not\subset CW_{\eps,K}^s(x)\qquad \mbox{ and }\qquad
CB_\delta(x)\not\subset CW_{\eps,K}^u(x)$$
\end{lema}
We are now in conditions to prove
\begin{teo}
If $X$ is a locally connected compact metric space, then for each
$\eps>0$ there exists $\delta>0$ such that
$\inf_{x\in X'}{\rm diam}\:CW^\sigma_\eps(x)\geq\delta$ for
$\sigma=s,u$.
\end{teo}
\setcounter{section}{2} \bp Let $0<\eps<\alpha$, and
choose $0<\delta<\eps$ so that, for each $x\in X$
\begin{equation}\label{delta.eps}
 d(x,y)\leq
\delta \quad\Rightarrow \quad d(f^n(x)f^n(y))< \eps\qquad \mbox{
for all } |n|\leq N_\eps\end{equation}
where $N_\eps$ is as in Lemma \ref{cw-uniform}. Let us denote
$CW_{\eps,n}^s(x)=cc(W_{\eps,n}^s(x),x)$, for each $n\geq 0$. We
shall see that $CW_{\eps,n}^s(x)\cap
\partial B_\delta(x)$ is not empty if $n\geq
\max\{N_\eps,K_{\eps,\delta}\}$. From this it follows that
$CW^s_\eps(x)=\lim_H CW_{\eps,n}^s(x)$ satisfies ${\rm diam}(CW_\eps^s(x))\geq \delta$.\newline\par
From Lemma \ref{estabilidad.finita} we have that $f^{-n}(CB_\delta(f^n(x)))$ is not included in
$CW_{\eps,n}^s(x)$ if $n\geq K_{\eps,\delta}$. Due to connectedness, there is a point
$y_0$ in $f^{-n}(CB_\delta(f^n(x)))$ belonging to $\partial
CW_{\eps,n}^s(x)$ which is included in $\partial W^s_{\eps,n}(x)$ because of local connectedness of $X$. This
implies
$$d(f^k(x),f^k(y_0))=\eps \qquad\quad\mbox{ for
some }k=0,\dots,n-1$$ Now, if $n\geq N_\eps$, then
\begin{itemize}
\item $k$ cannot be in $[n-N_\eps,n-1]$ because of (\ref{delta.eps}).
\item $k$ cannot be in $[N_\eps,n-N_\eps]$, for
otherwise there would exist $0\leq j\leq n$ such that
${\rm diam}(CW_\eps(x))>\alpha>2\eps$.
\end{itemize}
Therefore $k$ is mandatorily in $[0,N_\eps]$, resulting from
(\ref{delta.eps}) that $d(x,y_0)>\delta$.\ep
%______________________________________________________________________________________________
\section{Example}
%______________________________________________________________________________________________
\label{ejemplo.1} Here we show that $X$ being a continuum is not a
sufficient condition for an expansive $f$ to have local stable or
unstable sets of non trivial diameter for {\em all} points in
$X$.\par
\begin{figure}[t]\begin{center}
\psfrag{p}{$\tilde p$}\psfrag{a}{(a)}\psfrag{b}{(b)}
\includegraphics[width=12cm]{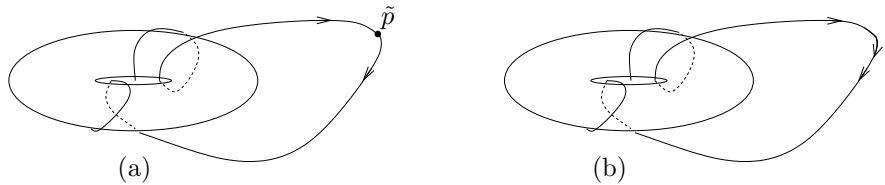}
\caption{\label{ejemplo.conexo}Examples}
\end{center}
\end{figure}
Take a linear Anosov automorphism $f$ on ${\Bbb T}^2$, and choose
a fixed point $p\in {\Bbb T}^2$. Let us also choose a stable and
an unstable separatrix of $p$, say $W^s_+(p)$ and $W^u_+(p)$. Lift
copies $\tilde p$, $\widetilde { W^s_+(p)}$, and
$\widetilde{W^u_+(p)}$ of $p$, $W^s_+(p)$ and $W^u_+(p)$
respectively, in ${\Bbb R}^3$, so that $\tilde{p}\in
\widetilde{W^s_+(p)}\cap\widetilde{W^u_+(p)}$ and
$\widetilde{W^\sigma_+(p)}$ be asymptotic to $W^\sigma_+(p)$,
$\sigma =s,u$ (see Figure \ref{ejemplo.conexo}.a). We can
introduce an expansive dynamic on this set extending the original
hyperbolic diffeomorphism. In this way we obtain an expansive
homeomorphism with stable points (in the future and the past),
none of which are fixed or periodic. This is essentially the
example of \cite{rr}.\par If $G$ is any small local
perturbation of $F$ at $\tilde p$, so that there are no fixed
points in $L=\widetilde{W^s_+(p)}\cup\widetilde{W^u_+(p)}$, and
$\omega(x)\cup\alpha(x)\subset {\Bbb T}^2$ for each point in $L$,
then $G$ is expansive, and has no stable points, while
$W^s_\eps(x)\cup W^u_\eps(x)=\{x\}$ for all $x$ in $L$ (see Figure
\ref{ejemplo.conexo}.b).
%_______________________________________________________________________


\begin{thebibliography}{RR}
%_______________________________________________________________________
\bibitem[H]{h}{\sc Hiraide K, }{Expansive homeomorphisms of compact
surfaces are Pseudo Anosov}, Osaka J. Math {\bf 27}, 1, (1990)
117-162.
\bibitem[K1]{ka}{\sc Kato H,} {Continuum-wise expansive homeomorphisms},
 Canad. J. Math. {\bf 45}, no.3, (1993) 576--598.
\bibitem[L1]{l1}{\sc Lewowicz J,} {Persistence in expansive systems},
Erg. Th. \& Dyn. Sys. {\bf 3} (1983) 567-578.
\bibitem[L2]{l2}{\sc Lewowicz J, } {Expansive homeomorphisms on surfaces},
Bol. Soc. Bras. Mat., Nova Ser. {\bf 20}, no. 1, (1989) 113-133.
\bibitem[P]{m}{\sc Paternain M, } {The principal loop-bundle and dynamical
systems}, C. R. Acad. Sci. Paris S\'{e}r. I Math. 329 (1999), no. 12,
1081-1085.
\bibitem[RR]{rr}{\sc Reddy W, Robertson L, }{Sources, sinks and saddles
for expansive homeomorphism with canonical
coordinates}, Rocky Mt. J. Math. {\bf 17} (1987) 673-681
\bibitem[V]{v} {\sc Vieitez J, } { Lyapunov Functions and Expansive
Diffeomorphisms on $3D$-Manifolds}, Ergodic Theory Dinam. Sys 22
(2002), no. 2, 601-632
\end{thebibliography}
\end{document}